\documentclass[11pt,a4paper,reqno]{amsart}
\usepackage{amsfonts}
\usepackage{graphicx}   
\usepackage{graphics}
\usepackage{epsfig}
\input{psfig.sty}
\usepackage{anysize}
\marginsize{2.5cm}{2.5cm}{2.5cm}{2.5cm}

\def\({\left(}
\def\){\right)}

\begin{document}
\bigskip
\bigskip
\title[]
{On one inverse spectral problem relatively domain}
\author[]{Y.S. Gasimov \dag}
\thanks{\dag Institute of Applied mathematics Baku State University, Z. Khalilov 23,
AZ1148 Baku, Azerbaijan, e-mail:\underline {gasimovyusif@bsu.az}}

\begin{abstract}
Different practical problems, espesially, problems of
hydrodynamics, elasticity theory, geophysics and aerodynamics can
be reduced to finding of an optimal shape. The investigation of
these problems is based on the study of depending of the
functiuonals on the domain, their first variation and gradient.

In the paper the inverse problem relatively domain is considered
for two-dimensional Schrodinger operator and operator $Lu = \Delta
^{2}u$ and the definition of $s - $functiuons is introduced. The
method is proposed for the determination of the domain by given
set of functions.



\end{abstract}
\maketitle

\bigskip

\textbf{1. Introduction}

One of the well studied classes of the inverse problems-is inverse spectral
problems. The papers dedicated to the investigation of these problems
traditionally focus on the construction of a function (potential) by given
spectral data (scattering data, normalizing numbers, eigenvalues…) and
obtaining necessary and sufficient conditions providing unequivocal
determination of the sought function. More datail review can be found in the
paper [1].

There exists a wide class of practical problems requiring determination of
the domain by some experimental data. For example, it is very important to
find the domain of the plate under across vibrations by the quantities,
which may be measured from distance [2]. There are different formulations of
the inverse problems relatively domain for the different ceses [3-5]. Note
that in differ from traditional inverse problems, the inverse problems
relatively domain have some special specifications. First, these problems
require to find no function, but domain. Second, choice of data (results of
measure) sufficient for determination of the domain is also enough difficult
problem.

In the paper we consider formulation and investigation of one inverse
problem relatively domain for the two dimensional Schrodinger operator that,
in particular, describes the vibration of membrane. In the end of the work
we put and solve the similar problem for the operator describing across
vibrations of the plate.

\bigskip

\textbf{2. Problem setting and preliminary results}

The object under investigation is the problem

\begin{equation}
\label{eq1}
 - \Delta u + q\left( {x} \right)u = \lambda u,\,\,\,\,x \in D,
\end{equation}

\begin{equation}
\label{eq2}
u\left( {x} \right) = 0,\,\,\,\,\,\,\,x \in S_{D} ,
\end{equation}

\noindent where $q(x)$ is differentiable non-negative function,
satisfying condition $t^{2}q\left( {xt} \right) = q\left( {x}
\right),\,\,\,t \in R,0 \notin D \subset R^{2}$-bounded convex
domain, $S_{D} \in C^{2} - $ its boundary, $\Delta $-Laplace
operator.

It is known [6,p.333] that under these conditions eigenfuntions
$u_{j}(x)$ of the problem (\ref{eq1}), (\ref{eq2}) belong to the
class $C^{2}\left( {D} \right) \cap C^{1}\left( {\bar {D}}
\right)$ and eigenvalues $\lambda _{j} $ are positive and may be
numbered as $\lambda _{1} \le \lambda _{2} \le ...$ considering
their multipilicity.

The set of all convex bounded domains $D \in R^{2}$ we denote by
$M$. Let

\[
K = \left\{ {D \in M:S_{D} \in \dot {C}^{2}} \right\},
\]

\noindent
where $\dot {C}^{2}$ is a class of the piece-wise twice continuous
differentiable functions.

\textbf{Definition1.} The functions

\begin{equation}
\label{eq3} J_{j} \left( {x,D} \right) = \frac{{\left| {\nabla
u_{j} \left( {x} \right)} \right|^{2}}}{{\lambda _{j}} },\,\,\,\,x
\in E^{2} ,\,\,\,j = 1,2, \ldots
\end{equation}

\noindent are said to be $s$- functions of the problem
(\ref{eq1}), (\ref{eq2}) in the domain $D$.

The problem is: To find a domain $D \in K$, such that

\begin{equation}
\label{eq4}
J_{j} \left( {x,D} \right) = s_{j} \left( {x} \right),\,\,\,\,x \in S_{D}
,\,\,\,j = 1,2, \ldots ,
\end{equation}

\noindent where $u_{j} \left( {x} \right)$ and $\lambda _{j} $ is
an eigenfunction and eigenvalue of the problem
(\ref{eq1})-(\ref{eq2}) in the domain $D$ correspondingly, $s_{j}
\left( {x} \right)$-given continuous functions defined\textbf{} on
$R^{2}$.

First of all we give the considerations, which led us to this formulation.

Investigation of the dependence of the eigenvalues of the operators on the
domain is an important problem, as the mechanical characteristics of some
systems indeed are eigenvalues of the corresponding operators, which may be
expressed by the functionals depending on the domain [7,8]. One of the
important steps for the investigation of the properties of these
characteristics is calculation of the variation of these functionals
relatively domain.

But to do this we need to define the space of domains, give a scalar product
and a definition for the domain variation in that space.

It was shown [9] the pairs $\left( {A,B} \right) \in M \times M$ form a
linear space with operations

\[
\left( {A,B} \right) + \left( {C,D} \right) = \left( {A + C,B + D}
\right),
\]

\[
 (A,B)=(C,D),\,\,\,\,   if \,\,\, \,\,   A + B = C + D,
\]

\[
\lambda \left( {A,B} \right) = \left( {\left| {\lambda}  \right|A,\left|
{\lambda}  \right|B} \right).
\]

Here $A + B$ is taken in the sence of Minkowsky, i.e.

\[
A + B = \left\{ {a + b:\,\,\,a \in A,\,\,\,b \in B} \right\}.
\]

The scalar product is introduced by formula

\[
(a,b) = \int\limits_{S_B }^{} {P_1^{} (\xi )P_2 (\xi )d\xi ,}
\]

\noindent
where

\[
a = \left( {A_{1} ,A_{2}}  \right),\,\,b = \left( {B_{1} ,B_{2}}  \right),
\]

\[
P_{1} \left( {x} \right) = P_{A_{1}}  \left( {x} \right) - P_{A_{2}}  \left(
{x} \right),
\]

\[
P_{2} \left( {x} \right) = P_{B_{1}}  \left( {x} \right) - P_{B_{2}}  \left(
{x} \right),
\]

 $S_{B} -$ unit sphere, $
P_D \left( x \right) = \mathop {\max }\limits_{l \in D} \left(
{x,l} \right),x \in E^2  - $ support function of the domain $D$.

The obtained space we define by $ML_{2} $.

For any fixed $D \in M$ eigenvaulue $\lambda _{j} $ of the problem
(\ref{eq1}), (\ref{eq2}) is defined as ([10], p.182)

\[
\lambda _{j} = infI\left( {u,D} \right),\,\,\left( {u,u_{p}}
\right) = 0,\,\,p = \overline {1,j - 1} ,
\]

\noindent
where

\[
I\left( {u,D} \right) = \frac{{\int\limits_{D} {\left[ {\left| {\nabla
u\left( {x} \right)} \right|^{2} + q\left( {x} \right)u^{2}\left( {x}
\right)} \right]dx}} }{{\int\limits_{D} {u^{2}\left( {x} \right)dx}} }.
\]

Thus we can consider $\lambda _{j} $ as a functional of $D \in K$
and define by $\lambda_{j}(D)$. The following formula is obtained
(see [9], p.98) for the first variation of the functional $\lambda
_{j} \left( {D} \right)$ in the space $ML_{2} $

\begin{equation}
\label{eq5}
\delta \lambda _{j} \left( {D} \right) = - \mathop {max}\limits_{u_{j}}
\int\limits_{S_{D}}  {\left| {\nabla u_{j} \left( {\xi}  \right)}
\right|^{2}\delta P_{D} \left( {n\left( {\xi}  \right)} \right)ds} ,
\end{equation}

\noindent where $\left| {\nabla u_{j} \left( {x} \right)}
\right|^{2} = \sum\limits_{i = 1}^{2} {\left( {\frac{{\partial
u\left( {x} \right)}}{{\partial x_{i}} }} \right)} ^{2}$,
$n(\xi)$- outside normal to $S_{D} $ in the point $\xi $, $max$ is
taken over all eigenfunctions corresponding to the eigenvalue
$\lambda _{j} $ in the case of its multipilicity.

Using (\ref{eq5}) the following formula may be obtained for the
eigenvalues of the problem (\ref{eq1}), (\ref{eq2}) in the domain
$D$

\begin{equation}
\label{eq6}
\lambda _{j} \left( {D} \right) = \frac{{1}}{{2}}\mathop {max}\limits_{u_{j}
} \int\limits_{S_{D}}  {\left| {\nabla u_{j} \left( {\xi}  \right)}
\right|^{2}P_{D} \left( {n\left( {\xi}  \right)} \right)ds} .
\end{equation}

Really, let's take $D_{0} \in K,\,\,D\left( {t} \right) = t \cdot D_{0}
,\,\,\,\,t > 0.$

By $u_{j} $ we define $j$-th eigenfunction of the problem
(\ref{eq1}), (\ref{eq2}) corresponding to the domain $D_{0}$. Then

\[
 - \Delta u_{j} \left( {x} \right) + q\left( {x} \right)u_{j} \left( {x}
\right) = \lambda _{j} u_{j} \left( {x} \right),\,\,x \in D_{0} .
\]

This relation may be written in the following equivalent from

\begin{equation}
\label{eq7}
 - \frac{{1}}{{t^{2}}}\Delta _{t} u_{j} \left( {\frac{{x}}{{t}}} \right) +
\frac{{1}}{{t^{2}}}q\left( {\frac{{x}}{{t}}} \right)u_{j} \left(
{\frac{{x}}{{t}}} \right) = \frac{{\lambda _{j} \left( {D_{0}}
\right)}}{{t^{2}}}u_{j} \left( {\frac{{x}}{{t}}} \right),
\quad
x \in D\left( {t} \right).
\end{equation}

Since the function
\[
 \tilde {u}_{j} \left( {x} \right) = u_{j} \left(
{\frac{{x}}{{t}}} \right),\,\,x \in D\left( {t} \right)
\]

\noindent
satisfies the relation

\begin{equation}
\label{eq8} \Delta \tilde u_j (x) = \frac{1}{{t^2 }}\Delta u_j
\left( {\frac{x}{t}} \right)
\end{equation}

\noindent
from the condition $t^{2}q\left( {tx} \right) = q\left( {x} \right)$and (8)
one may get

\[
 - \Delta \tilde {u}_{j} \left( {x} \right) + q\left( {x} \right)\tilde
{u}_{j} \left( {x} \right) = \frac{{\lambda _{j} \left( {D_{0}}
\right)}}{{t^{2}}}\tilde {u}_{j} \left( {x} \right),\,\,x \in D\left( {t}
\right).
\]

It shows that $\Delta \tilde {u}_{j} \left( {x} \right)$ is an
eigefunction, and $\lambda _{} \left( {t} \right) = \frac{{\lambda
_{j} \left( {D_{0}} \right)}}{{t^{2}}}$ eigenvalue for the problem
(\ref{eq1}), (\ref{eq2}) in the domain $D(t)$. Then using
(\ref{eq5}) we can write

\begin{equation}
\label{eq8}
\begin{array}{l}
 \,\,\,\,\,\,\,\,\,\,\,\,\,\,\,\,\,\,\,\,\,\lambda _{j} \left( {t + \Delta
t} \right) - \lambda _{j} \left( {t} \right) = \lambda _{j} \left( {D\left(
{t + \Delta t} \right)} \right) - \lambda _{j} \left( {D\left( {t} \right)}
\right) = \\
 = \int\limits_{S_{D\left( {t} \right)}}  {\left| {\nabla u\left( {\xi}
\right)} \right|^{2}\left[ {P_{D\left( {t + \Delta t} \right)}
\left( {n\left( {\xi}  \right)} \right) - P_{D\left( {t} \right)}
\left( {n\left( {\xi}  \right)} \right)} \right]} ds + o\left(
{\Delta t} \right),\,\,\,\xi
\in S_{D\left( {t} \right)} \\
 \end{array}.
\end{equation}

If support function $P_{D\left( {t} \right)} \left( {x} \right)$
of the domain $D\left( {t} \right)$ is differentable relatively
$t,$ then dividing both sides of (\ref{eq8}) by $t$ we obtain

\begin{equation}
\label{eq9}
\lambda^\prime _j (t) =  - \mathop {\max
}\limits_{u_j}
 \int\limits_{S_{D\left( {t} \right)}}  {\left| {\nabla u_{j}\left( {x}
\right)} \right|^{2}{P}'_{D(t)} (n(\xi)) ds},
\end{equation}

\noindent where ${P}'_{D\left( {t} \right)} \left( {x} \right) =
\frac{{\partial }}{{\partial t}}P_{D\left( {t} \right)} \left( {x}
\right)$.

Considering this we have

\[
 - 2\frac{{\lambda _{j} \left( {D_{0}}  \right)}}{{t^{3}}} = -
\frac{{1}}{{t^{2}}}\mathop {max}\limits_{u_{j} \left( {x} \right)}
\int\limits_{S_{D}}  {\left| {\nabla u_{j} \left( {\frac{{\xi} }{{t}}}
\right)} \right|^{2}P_{D_{0}}  \left( {n\left( {\xi}  \right)}
\right)dS,\,\,\,\xi \in S_{D}}  .
\]

Taking $t = 1$ from this we get (\ref{eq6}).

As we see from (\ref{eq6}) the boundary values of the function
$\left| {\nabla u_j \left( x \right)} \right|^2 $ unquiovoccally
define eigenvalue $\lambda _{j} $.

From (\ref{eq6}) taking into account (\ref{eq4}) we obtain

\begin{equation}
\label{eq10}
\int\limits_{S_{D}}  {s_{j} \left( {\xi}  \right)P_{D} \left( {n\left( {\xi
} \right)} \right)ds = 2,\,\,\,\,j = 1,2, \ldots}  \quad .
\end{equation}

This is the basic relation for the solving of the considered
problem.

\textbf{Note.} As we take $s$- functions as a given data, let's consider
them for some concrete cases. For one dimensional case

\begin{equation}
\label{eq11}
{u}'' + q\left( {x} \right)y = \lambda u,
\end{equation}

\begin{equation}
\label{eq12} u(a)=u(b)=0,
\end{equation}

\noindent
where $q\left( {x} \right) = \frac{{c}}{{x^{2}\,}},\,\,c \ge 0,\,\,0 \notin
\left( {a,b} \right) \subset R$, $s$-functions defined by (\ref{eq3}) indeed are

\[
\begin{array}{l}
 \frac{{u_{jx}^{2} \left( {a} \right)}}{{\lambda _{j}} } = J_{j} \left( {a}
\right), \\
 \frac{{u_{jy}^{2} \left( {b} \right)}}{{\lambda _{j}} } = J_{j} \left( {b}
\right). \\
 \end{array}
\]

Thus the expression (8) takes a form

\begin{equation}
\label{eq13} J_{j} \left( {b} \right) \cdot b - J_{j} \left( {a}
\right) \cdot a = 2,j = 1,2,...
\end{equation}

Let's take $a=0$, i.e. consider the problem (\ref{eq11}),
(\ref{eq12}) in the interval $\left( {0,b} \right)$. For this case
from (\ref{eq8}) we get the following

\textbf{Consequence}. All $s$-functions of the problem (\ref{eq11}), (\ref{eq12}) satisfy to
the condition

\begin{equation}
\label{eq14} J_{j} \left( {b} \right) = \frac{{2}}{{b}},\,\,j =
1,2,... \quad .
\end{equation}

This formula allows to solve the inverse problem: Let the set of
functions $s_{j} \left( {x} \right),\,\,j = 1,2, \ldots $ is
given. In this case the problem of finding of the domain
satisfying (\ref{eq4}) is reduced to determination of the point
$b$, which may be done using (\ref{eq14}).

As noted in consequence all $s$-functions satisfy to the condition
(\ref{eq14}), which is equivalent to the one condition. This
condition is sufficient for finding of the point $b$. Really as
one may get from (\ref{eq14})

\[
b = \frac{{2}}{{J_{j} \left( {b} \right)}}.
\]

Similarly, if $b = 0$, then we have

\[
a =  - \frac{2}{{J_j (a)}}.
\]

Note that if $J_{j} \left( {x} \right) \equiv c_{j} $, $x \in
S_{D} $, $c_{j} = const,j = 1,2, \ldots $ then as follows from
(\ref{eq14}) they all are equal to each other for all $j =
1,2,...$ .

In two dimensional case from (\ref{eq10}) is obtained that if the
functions $J_{j}(x,D)$ are constant, then

 $J_{j} \left( {x,D} \right) \equiv \frac{{1}}{{mesD}},\,\,\,\,j = 1,2, \ldots
$ (see [9]).

Now we prove the lemma that will be used later on.

\textbf{Lemma1.} Let $f\left( {x} \right)$ be continuous function defined on
the unit shere$S_{B} $. Then for any $D_{1} ,\,\,D_{2} \in K$

\begin{equation}
\label{eq15} \int\limits_{S_{D_1  + D_2 } } {f(n} (\xi ))ds =
\int\limits_{S_{D_1 } } {f(n(\xi ))ds + } \int\limits_{S_{D_2 } }
{f(n(\xi ))ds} ,
\end{equation}

\noindent
where $D_{1} + D_{2} $ is taken in the sence of Minkowsky i.e.

\[
D_{1} + D{}_{2}^{} = \left\{ {x:x = x_{1} + x_{2} ,\,x{}_{1}^{} \in D_{1}
,\,x_{2} \in D_{2}}  \right\}.
\]

\textbf{Proof.} It is known [11], that $f\left( {x} \right)$ may be
continuously, positive-homogeneously extended over all the space and
presented as a limit of the difference of two convex functions

\begin{equation}
\label{eq16}
f\left( {x} \right) = \mathop {lim}\limits_{n \to \infty}  \left[ {g_{n}
\left( {x} \right) - h_{n} \left( {x} \right)} \right].
\end{equation}

First consider

\begin{equation}
\label{eq17} f(x)=g(x)-h(x),
\end{equation}

\noindent
where $g\left( {x} \right),\,h\left( {x} \right)$ are convex,
positively-homogeneous functions.

As is known [12] for any continuous, convex,
positively-homogeneous function $P\left( {x} \right)$ there exists
the only convex bounded set $D$ such, that $P\left( {x} \right)$
is a support function of $D$, i.e. $P\left( {x} \right) = P_{D}
\left( {x} \right)$. The opposite statement also is true.

It is also known that $D$ is found as subdifferential of its support
function in the point $x = 0$

\[
D = \partial P\left( {0} \right) = \left\{ {l \in E^{n}:P\left( {x} \right)
\ge \left( {l,x} \right),\,\forall x \in R^{n}} \right\}.
\]

So, there exist the domains $G$ and $H$, such that

\begin{equation}
\label{eq18}
g\left( {x} \right) = P_{G} \left( {x} \right),\,\,\,h\left( {x} \right) =
P_{H} \left( {x} \right).
\end{equation}

Considering (\ref{eq17}), (\ref{eq18}) we get

\begin{equation}
\label{eq19}
\begin{array}{l}
 \int\limits_{S_{D_{1} + D_{2}} }  {f\left( {n\left( {\xi}  \right)}
\right)ds =}  \int\limits_{S_{D_{1} + D_{2}} }  {\left[ {g\left( {n\left(
{\xi}  \right)} \right) - h\left( {n\left( {\xi}  \right)} \right)ds}
\right] =}  \\
 \\
 \,\,\,\,\,\,\,\, = \int\limits_{S_{D_{1} + D_{2}} }  {P_{G} \left( {n\left(
{\xi}  \right)} \right)ds -}  \int\limits_{S_{D_{1} + D_{2}} }  {P_{H}
\left( {n\left( {\xi}  \right)} \right)ds.} \\
 \end{array}
\end{equation}

As for any$D_{1} ,D_{2} \in K$ the following relation is valid [9]

\begin{equation}
\label{eq20} \int\limits_{S_{D_1 } } {P_{D_2 } \left( {n\left( \xi
\right)} \right)ds = } \int\limits_{S_{D_2 } } {P_{D_1 } \left(
{n\left( \xi  \right)} \right)ds} ,
\end{equation}

\noindent
from (\ref{eq19}) one may obtain

\[
\int\limits_{S_{D_{1} + D_{2}} }  {f\left( {n\left( {\xi}  \right)}
\right)ds =}  \int\limits_{S_{G}}  {P_{D_{1} + D_{2}}  \left( {n\left( {\xi
} \right)} \right)ds} - \int\limits_{S_{H}}  {P_{D_{1} + D_{2}}  \left(
{n\left( {\xi}  \right)} \right)ds} .
\]

As $P_{D_{1} + D_{2}}  \left( {x} \right) = P_{D_{1}}  \left( {x} \right) +
P_{D_{2}}  \left( {x} \right)$ [12], applying (\ref{eq20}) again we get (\ref{eq15}). Lemma
is proved.

\bigskip

\textbf{3. Main results}

Now we investigate the main problem of the work-construction of $D$ by given
set of functions $s_{j} \left( {x} \right),\,\,\,j = 1,2, \ldots $.

Let $B \subset E^2$ be unit ball with the center at the origin and
$S_{B} $- its boundary. By $\varphi _{k} \left( {x} \right),\,\,k
= 1,2, \ldots $ we denote some basis in $C\left( {S_{B}}
\right)$-space of continuous in $S_{B} $ functions. These
functions may be continuously, positive-homogeneously extended to
$B$. It may be done as:

\[
\tilde \varphi _k (x) = \left\{ \begin{array}{l}
 \varphi _k \left( {\frac{x}{{\left\| x \right\|}}} \right) \cdot \left\| x \right\|,x \in B,x \ne 0, \\
 0,x = 0. \\
 \end{array} \right.
\]

One may check, that these functions are continuous and satisfy to the
positive homogeneity condition

\[
\tilde {\varphi} _{k} \left( {\alpha x} \right) = \alpha \tilde {\varphi
}_{k} \left( {x} \right),\,\,\,\alpha > 0.
\]

Without loss of generality we can denote $\tilde {\varphi} _{k}
\left( {x} \right)$ by $\varphi _{k} \left( {x} \right)$.

Thus we obtain the set of continuous, positive-homogeneous
functions defined on $B$.

As we noted above each positive-homogeneous, continuous function $\varphi
_{j} \left( {x} \right)$may be presented in the form

\begin{equation}
\label{eq21}
\varphi _{k} \left( {x} \right) = \mathop {lim}\limits_{n \to \infty}
\left[ {g_{n}^{k} \left( {x} \right) - h_{n}^{k} \left( {x} \right)}
\right]
\end{equation}

\noindent and there exist satisfying above mentioned properties
domains $G_{n}^{k} $ and $H_{n}^{k} $ such, that

\[
g_{n}^{k} \left( {x} \right) = P_{G_{n}^{k}}  \left( {x} \right),
\quad
h_{n}^{k} \left( {x} \right) = P_{H_{n}^{k}}  \left( {x} \right).
\]

These domains we call basic domains. Substituting these into (\ref{eq21}) we get

\begin{equation}
\label{eq22}
\varphi _{k} \left( {x} \right) = \mathop {lim}\limits_{n \to \infty}
\left[ {P_{G_{n}^{k}} ^{} \left( {x} \right) - P_{H_{n}^{k}} ^{} \left( {x}
\right)} \right].
\end{equation}

First we consider

\begin{equation}
\label{eq23}
\varphi \left( {x} \right) = P_{G^{k}} \left( {x} \right) - P_{H^{k}} \left(
{x} \right),
\end{equation}

\noindent where $G^{k}$and $H^{k}$ are closed, bounded convex
domains.

As $n\left( {x} \right) \in S_{B} $, for any $x \in S_{D} $, we
can decompose $P_{D} \left( {x} \right)$,$x\in S_{B}$ by basic
functions $\varphi _{k} \left( {x} \right)$

\begin{equation}
\label{eq24}
P_{D} \left( {x} \right) = \sum\limits_{k = 1}^{\infty}  {\alpha _{k}
\varphi _{k} \left( {x} \right),\,\,\,x \in S_{B}}  ,\,\,\,\alpha _{k} \in
R,\,\,\,k = 1,2, \ldots .
\end{equation}

Thus, to determine $P_{D} \left( {x} \right)$ we have to find the
coefficients $\alpha _{k} ,\,\,k = 1,2, \ldots $ .

\textbf{Theorem1.} Let the set of functions $s_{j}(x), j=1,2,...$
is given. Then the coefficients $\alpha _{k} ,\,\,\,\,k = 1,2,
\ldots $ of the support function of sought domain $D$ for which
(\ref{eq4}) is true, satisfy the equation

\begin{equation}
\label{eq25}
\sum\limits_{k,m = 1}^{\infty}  {A_{k,m} \left( {j} \right)\alpha _{k}
\alpha _{m} = 2,\,\,\,\,\,j = 1,2, \ldots}
\end{equation}

\noindent
with coefficients

\begin{equation}
\label{eq26}
\begin{array}{l}
 A_{k,m} \left( {j} \right) = \int\limits_{S_{G^{k}}}  {s_{j} \left( {x}
\right)\left[ {P_{G^{m}} \left( {n\left( {x} \right)} \right) - P_{H^{m}}
\left( {n\left( {x} \right)} \right)} \right]ds} - \\
 - \int\limits_{S_{H^{k}}}  {s_{j} \left( {x} \right)\left[ {P_{G^{m}}
\left( {n\left( {x} \right)} \right) - P_{H^{m}} \left( {n\left( {x}
\right)} \right)} \right]ds} . \\
 \end{array}
\end{equation}

\textbf{Proof.} Considering (\ref{eq23}) from (\ref{eq24}) one may get

\begin{equation}
\label{eq27} P_D (x) = \sum\limits_{k = 1}^\infty  {\alpha _k
(P_{G^k } (x) - P_{H^k } (x))} ,x \in S_B .
\end{equation}

The set of all indexes for which $\alpha _{k} \ge 0  \left(
{\alpha _{k} < 0} \right)$ denote by $I^{ + \,\,}\left( {I^{ -} }
\right)$.

Then the relation (\ref{eq27}) may be written as

\begin{equation}
\label{eq28}
\begin{array}{l}
 P_{D} \left( {x} \right) - \sum\limits_{k \in I^{ -} }^{} {\alpha _{k}
P_{G^{k}} \left( {x} \right) + \sum\limits_{k \in I^{ +} } {\alpha _{k}
P_{H^{k}}}  \left( {x} \right) =}  \\
 \,\, = \sum\limits_{k \in I^{ +} }^{} {\alpha _{k} P_{G^{k}} \left( {x}
\right) - \sum\limits_{k \in I^{ -} } {\alpha _{k} P_{H^{k}} \left( {x}
\right),\,\,\,\,\,}}  x \in S_{B} . \\
 \end{array}
\end{equation}

From last taking into account the properties of support functions [12] we
obtain

\[
D - \sum\limits_{k \in I^ -  }^{} {\alpha _k G^k  + \sum\limits_{k
\in I^ +  } {\alpha _k H^k  = } } \sum\limits_{k \in I^ +  }^{}
{\alpha _k G^k  - \sum\limits_{k \in I^ -  } {\alpha _k H^k } }.
\]

The use of (\ref{eq28}) and the lemma gives

\[
\begin{array}{l}
 \int\limits_{S_D } {s_j } (\xi )P_D (n(\xi )d\xi  + \int\limits_{\sum\limits_{k \in I^ -  } {( - \alpha _k )S_{G^k } } } {s_j (\xi )P_D (n(\xi ))d\xi  + }  \\
  \\
  + \int\limits_{\sum\limits_{k \in I^ +  } {\alpha _k } S_{H^k } } {s_j } (\xi )P_D (n(\xi )d\xi  = \int\limits_{\sum\limits_{k \in I^ +  } {\alpha _k S_{G^k } } } {s_j (\xi )P_D (n(\xi ))d\xi }  +  \\
  \\
  + \int\limits_{\sum\limits_{k \in I^ -  } {( - \alpha _k )S_{H^k } } } {s_j } (\xi ))P_D (n(\xi ))d\xi . \\
 \end{array}
\]

From this considering (\ref{eq10}) we have

\[
\begin{array}{l}
 \int\limits_{S_D } {s_j } (\xi )P_D (n(\xi )d\xi  = \sum\limits_{k = 1}^\infty  {\alpha _k \left[ {\int\limits_{S_{G^k } } {s_j } } \right.} (\xi )P_D (n(\xi ))d\xi  -  \\
 \left. { - \int\limits_{S_{H^k } } {s_j } (\xi )P_D (n(\xi )d\xi } \right] = 2. \\
 \end{array}
\]

Substituting here (\ref{eq24}) one may get (\ref{eq25}) with coefficients (\ref{eq26}). Theorem is
proved.

We assumed that the considered problem has a solution in general case. For
interesting cases as the functions $s_{j} \left( {x} \right),\,\,j = 1,2,
\ldots $ are defined as experimental data, this problem always has a
solution. The function $P_{D} \left( {x} \right)$ is constructed by the help
of the solution of (\ref{eq25}), using (\ref{eq24}).

As we noted above domain $D$ is unequivocally defined by its
support function $P_{D} \left( {x} \right)$. Suppose that
(\ref{eq25}) has the only solution providing convexity of the
support function of $D$.

Let's show that the expressions $\frac{{\left| {\nabla u_{j}
\left( {\xi} \right)} \right|^{2}}}{{\lambda _{j}} },\,\,\,\,\,j =
1,2,...$ for the problem (\ref{eq2}), (\ref{eq3}) in the
constructed by the help of this solution, using (\ref{eq24})
domain $D$ indeed are $s$-functions. Really, if $D$ is a domain in
which the problem (\ref{eq2}), (\ref{eq3}) has given by formula
(\ref{eq24}) $s$-functions then decomposition $\overline {D} $ by
formulae (\ref{eq24}) and making above done transformations we get
the equation (\ref{eq25}) with the same coefficients. From the
assumption that this equation has the only solution, it follows
$\overline {D} = D$.

If (\ref{eq18}) has more than one solution then the sought domain
is among the ones, constructed by (\ref{eq17}) using these
solutions, providing convexity of $P(x)$.

This algorithm is constructed considering (\ref{eq17}). In general case when
$\varphi \left( {x} \right)$ has form (\ref{eq21}) $A_{k,m} \left( {j} \right)$
turns to

\[
\begin{array}{l}
 A_{k,m} \left( {j} \right) = \mathop {lim}\limits_{n \to \infty}  \left[
{\int\limits_{S_{G_{n} ^{k}}}  {s_{j} \left( {x} \right)\left[ {P_{G_{n}
^{m}} \left( {n\left( {x} \right)} \right) - P_{H_{n} ^{m}} \left( {n\left(
{x} \right)} \right)} \right]ds} -}  \right. \\
 \left. { - \int\limits_{S_{H_{n} ^{k}}}  {s_{j} \left( {x} \right)\left[
{P_{G_{n} ^{m}} \left( {n\left( {x} \right)} \right) - P_{H_{n} ^{m}} \left(
{n\left( {x} \right)} \right)} \right]ds}}  \right]. \\
 \end{array}
\]

Now let's consider the across vibrations of the plate.

Let $D \in R^{2}$ be a domain of the plate with boundary $ S_D \in
C^2 $.

It is known [2] that the function $\omega \,\left( {x_{1} x_{2} ,\,t}
\right)$ describing across vibrations of the plate satisfies equation

\begin{equation}
\label{eq29}
\omega _{x_{1} x_{1} x_{1} x_{1}}  + 2\omega _{x_{1} x_{1} x_{2} x_{2}}  +
\omega _{x_{2} x_{2} x_{2} x_{2}}  + \omega _{tt} = 0 \quad .
\end{equation}

Assuming the process stabilized the solution - eigen-vibration is sought as

\[
\omega \left( {x_{1} ,x_{2} ,t} \right) = u\left( {x_{1} ,x_{2}}
\right)cos\lambda t,
\]

\noindent where $\lambda$-is an eigen-frequency.

Substituting this to (\ref{eq29}) we get

\begin{equation}
\label{eq30}
\Delta ^{2}u = \lambda u,\,\,\,\,\,x \in D,\,\,\,
\end{equation}

\noindent
where $\Delta ^{2} = \Delta \Delta $.

For different cases different boundary conditions may be considered. The
object under investigation is the freezed plate with boundary conditions

\begin{equation}
\label{eq31}
u = 0,\,\,\,\,\frac{{\partial u}}{{\partial n}} = 0,\,\,\,\,\,x \in S_{D} .
\end{equation}

Let

\[
K = \left\{ {D \in M:S_{D} \in \dot {C}^{2}} \right\},
\]

\noindent where $ \dot C^2 $ is a class of the piece-wise twice
continuous differentiable functions.

\textbf{Definition2.} The functions $J_{j} \left( {x,D} \right) =
\frac{{\left| {\Delta u_{j} \left( {x} \right)}
\right|^{2}}}{{\lambda _{j} }},\,\,\,\,x \in E^{2} ,\,\,\,j = 1,2,
\ldots $ are called $s$- functions of the problem (\ref{eq30}),
(\ref{eq31}) in the domain $D$.

The problem is: To find a domain $D \in K$, such that

\begin{equation}
\label{eq32} J_j \left( {x,D} \right) = s_j (x),x \in S_D ,j =
1,2, \ldots,
\end{equation}

\noindent
where $u_{j} \left( {x} \right)$ and $\lambda _{j} $ are eigen-vibration and
eigen-frequency of the problem (\ref{eq30})-(\ref{eq31}) in the domain D correspondingly,
$s_{j} \left( {x} \right),\,\,j = 1,2, \ldots $-given continuous functions
defined\textbf{} on $R^{2}$.

In [9] the following formula is obtained for the eigen-frequency of the
freezed plate under across vibrations

\begin{equation}
\label{eq33} \lambda _j  = \frac{1}{4}\mathop {\max }\limits_{u_j
} \int\limits_{S_D } {\left| {\Delta u_j (\xi )} \right|^2 P_D
(n(\xi ))ds,}
\end{equation}

\noindent where $P_{D} \left( {x} \right) = \mathop
{max}\limits_{l \in D} \left( {l,x} \right),\,\,\,\,x \in E^{n}$
is a support function of $D$, and $max$ is taken over all
eigen-vibrations $u_{j} $ corresponding to eigen-frequency
$\lambda_{j}$ in the case of its multiplicity. (As we see from
(\ref{eq33}), the boundary values of the function $\left| {\Delta
u_{j} \left( {x} \right)} \right|^{2}$ unequivocally define
$\lambda _{j} $). From (\ref{eq33}) considering (\ref{eq32}) we
get

\[
\int\limits_{S_{D}}  {s_{j}}  \left( {\xi}  \right)P_{D} \left( {n\left( {x}
\right)} \right)ds = 4,\,\,\,j = 1,2,... \quad .
\]

Caring out above done considerations the following theorem is proved for
considered problem.

\textbf{Theorem2.} Let the set of functions $s_{j} \left( {x}
\right),\,\,j = 1,2, \ldots $ is given. Then the coefficients of
the support function of sought domain $D$ of the plate for which
(\ref{eq32}) is true, satisfy the equation

\[
\sum\limits_{k,m = 1}^{\infty}  {A_{k,m} \left( {j} \right)\alpha _{k}
\alpha _{m} = 4,\,\,\,\,\,j = 1,2,...} ,
\]

\noindent
with coefficients

\[
\begin{array}{l}
 A_{k,m} \left( {j} \right) = \mathop {lim}\limits_{n \to \infty}  \left[
{\int\limits_{S_{G_{n} ^{k}}}  {s_{j} \left( {x} \right)\left[ {P_{G_{n}
^{m}} \left( {n\left( {x} \right)} \right) - P_{H_{n} ^{m}} \left( {n\left(
{x} \right)} \right)} \right]ds} -}  \right. \\
 \left. { - \int\limits_{S_{H_{n} ^{k}}}  {s_{j} \left( {x} \right)\left[
{P_{G_{n} ^{m}} \left( {n\left( {x} \right)} \right) - P_{H_{n} ^{m}} \left(
{n\left( {x} \right)} \right)} \right]ds}}  \right]. \\
 \end{array}
\]

\bigskip

References:

1. Aktosun T., Weder R. Inverse spectral - scatterring problem
with two sets of discrete spectra for radial Schrodinger equation.
Preprint, 2004, 47 pp.

2. S.H. Gould. Variational methods for eigenvalue problems. Univ. of Toronto
Press. London: Oxford Univ. Press, 1996.

3. F Pesaint and J.-P. Zolesio. Derivees par rapport au domaine des valeurs
propres du Laplacien C.R. Acad.Sci.Ser.1 (1995), v.321, ¹10, p.1310-1337.

4. J.Elschner, G.Schimdt, M. Yamamoto. Au inverse problem in periodic
diffractive optics: global uniqueness with a single wave number. Inverse
problems, 2003, 19, p.779-787.

5. R. Dziri, J.-P. Zolesio. Shape derivative with Lipschits coindiuous
coefficients. Boll. Unione. Math. Ital. B (1996), v.10, N3, p.569-594.

6. V.S. Vladimirov. Equations of the Mathematical Physics. Nauka, Moscow,
1988 (in Russian).

7. Y.S. Gasimov, A.A. Niftiyev. On a minimization of the eigenvalues of
Schrodinger operator over domains. Reports of the Russian Academy of
Scienses, 2001, v.380, N3, p.305-307 (in Russian).

8. Y.S. Gasimov. On some properties of the eigenvalues when the domain
varies. Mathematical Physics, Analyses, Geometry, 2003, v.10, N2, p.249-255
(in Russian).

9. A.A. Niftiyev, Y.S. Gasimov. Control by boundaries and eigenvalue
problems with variable domain. Publishing House of Baku St. Univ., 2004 (in
Russian).

10. V.P. Mikhailov. Partial Differential Equations. Nauka, Moscow, 1976 (in
Russian).

11. D.M.Burago, V.A. Zalgamer. Geomethric inequalities. Nauka, Moskow,
1981(in Russian).

12. V.F. Demyanov, A.M. Rubinov. Basises of non-smooth analyses and
quazidifferential calculas. Nauka, Moscow , 1990.

\end{document}